\documentclass[runningheads]{llncs}

\usepackage[breaklinks, colorlinks, linkcolor=blue, citecolor=blue, urlcolor=black]{hyperref}
\usepackage{color}

\usepackage{graphicx}
\usepackage{tabularx}
\usepackage{tabularcalc}
\usepackage{booktabs}
\usepackage{multirow}
\usepackage{algorithm}
\usepackage[noend]{algpseudocode}
\usepackage{tikz}
\usepackage{forest}
\usepackage{pgfplots}
\usepackage{xltabular}
\usetikzlibrary{arrows.meta}
\usepackage{amsmath}
\usepackage{amssymb}
\usepackage{xspace}
\usepackage{xcolor}
\usepackage{varioref}
\usepackage{aligned-overset}
\usepackage[shortlabels,inline]{enumitem}
\usepackage{boxedminipage}
\usepackage{makecell}
\usepackage{underscore}
\usepackage{rotating}
\usepackage{cleveref}
\usepackage{hyperref}
\usepackage{pgfplots}
\usepackage{booktabs}
\usepackage{todonotes}
\usetikzlibrary{positioning}
\pgfplotsset{compat=1.18}

\newcommand{\RNum}[1]{\uppercase\expandafter{\romannumeral #1\relax}}
\newcommand{\solver}[1]{\textsc{#1}\xspace}

\newcommand{\chap}{\solver{CHAP}}

\newcommand{\myorcidlink}[1]{\,\href{https://orcid.org/#1}{\raisebox{-0.45ex}{\includegraphics[width=1.8ex]{orcid}}}}

\makeatletter
\def\orcidID#1{\href{http://orcid.org/#1}{\protect\raisebox{-1.25pt}{\protect\includegraphics{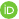}}}}
\makeatother

\begin{document}

\title{CHAP: A Hybrid GPU-CPU Heuristic for MIP}

\author{
    Gennesaret Kharistio Tjusila\inst{1} \orcidID{0009-0009-7031-9993} \and
    Alexander Hoen\inst{1,2}\orcidID{0000-0003-1065-1651} \and  
    Nils-Christian Kempke\inst{1}\orcidID{0000-0003-4492-9818} \and
    Gioni Mexi\inst{1}\orcidID{0000-0003-0964-9802}\and
    Timo Berthold\inst{1,3,4}\orcidID{0000-0002-6320-8154} \and
    Ambros Gleixner\inst{1,2}\orcidID{0000-0003-0391-5903} \and
    Thorsten Koch\inst{1,3}\orcidID{0000-0002-1967-0077} \and
    Sebastian Pokutta\inst{1,3}\orcidID{0000-0001-7365-3000}
}

\authorrunning{Tjusila et al.}

\institute{
    Zuse Institute Berlin, Takustraße 7, 14195 Berlin, Germany 
    \email{\{kempke,koch,mexi,tjusila,pokutta\}@zib.de}\and
    HTW Berlin, Berlin, Germany
    \email{\{gleixner,hoen\}@htw-berlin.de}\and
    Technische Universität Berlin, Berlin, Germany \and
    Fair Isaac Deutschland GmbH, Germany
    \email{timoberthold@fico.com}
    }

%  % place this right after \maketitle
% \begingroup
% \renewcommand\thefootnote{\dag}
% \footnotetext{These authors contributed equally to this work.}
% \endgroup

\maketitle             
\begin{abstract}
We present \chap (Coordinating Heuristics Across Plat\-forms) a GPU-CPU-hybrid primal heuristic framework for mixed-integer programming. 
\chap adopts a portfolio approach where it coordinates a set of primal heuristics, including Local Search, Fix-and-Propagate, and Feasibility Pump, via a shared solution pool. 
The solution pool is used to exchange feasible incumbent solutions, LP solutions, along with promising infeasible solution candidates, enabling a more comprehensive exploration of the solution space.
On the GPU side, we implement a native tabu search featuring a novel best-shift algorithm built on sort, scan, and reduce primitives, along with specialized kernel designs. We additionally leverage cuPDLPx as an approximate LP solver.
On the CPU side, we employ various Fix-and-Propagate strategies, guided by information from the solution pool, complemented by a CPU-based tabu search and a Feasibility Pump.
All components operate collaboratively, iteratively repairing and improving candidate solutions maintained in the pool. We evaluate our framework on the 50-instances benchmark from the 2026 Land--Doig MIP Competition under competition constraints, including a five-minute time limit. 
In these settings, \chap finds solutions to 47 instances outperforming both Gurobi (44) in default mode and NVIDIA cuOpt (43) in heuristics-only mode. 
The results demonstrate that coordinated cross-platform portfolios offer a promising direction for the integration of GPU heuristics into modern high-performance MIP solvers. 
The code will be made available on GitHub.
%
%\keywords{First keyword  \and Second keyword \and Another keyword.}
\end{abstract}

\section{Introduction}
\label{sec:introduction}

%\TBil{AG TB SP Considerations: Verteilungsproblemen, Autotuning, andere Hardware-Überlegungen (TK)}

Mixed-integer linear programming (MIP) is one of the central modeling and algorithmic paradigms in mathematical optimization.
At the same time, the computational landscape on which MIP solving methods operate is changing rapidly.
GPUs have moved beyond their traditional role as accelerators for isolated dense kernels. Recent progress on GPU-enabled first-order methods for large-scale linear optimization~\cite{applegate2021practical,applegate2025pdlp,lu2025geometry} shows that they can support substantial parts of the optimization pipeline.
Modern compute nodes are also heterogeneous within and across device classes, combining CPUs with different core types and vector units as well as GPUs with widely varying throughput and memory characteristics.
As a result, designing high-performance primal heuristics is no longer only a question of choosing the right search strategy; it is also a question of deciding which work should be carried out where, and how information should be exchanged without introducing excessive synchronization overhead.

A natural next step is therefore not merely to port an existing MIP heuristic to the GPU, but to design cross-platform heuristic frameworks in which CPU and GPU components reinforce one another.
This question is timely in light of recent GPU-based primal heuristics for MIP~\cite{kempke2026FP,corduk25GpuAcceleratedPrimalHeuristics} and is directly reflected in the 2026 Land--Doig MIP Competition~\cite{mipcompetition2026}, which focuses on GPU-accelerated primal heuristics for general MIP.

In this work, we present \chap, a primal-heuristic framework developed for this competition.
\chap adopts a portfolio view of heuristic search: complementary CPU and GPU workers exchange linear programming (LP) iterates, feasible incumbents, and promising infeasible solution candidates through a shared solution pool.
With a slight abuse of notation, we will refer to such infeasible candidate vectors as \emph{infeasible solutions} for the remainder of the paper.
At the core of this interaction are a GPU-native tabu search and a fix-and-propagate-and-repair~\cite{Salvagnin2024FPR} framework guided by streamed LP information.
This coordination is essential to the method: the GPU side provides fast large-neighborhood improvement, while the CPU side contributes repair, warmstarting, and intensification mechanisms that exploit the information produced across the portfolio.

Our contributions span both algorithm design and systems engineering.
Algorithmically, we combine a feasibility-jump-style routine~\cite{Luteberget2023FeasJump} with LocalMIP-inspired search~\cite{Lin25LocalMip}, introduce guided repair and warmstarting mechanisms, and study improved neighborhood-evaluation strategies.
On the systems side, we develop specialized GPU kernels, refine move-evaluation routines, and incorporate execution- and memory-level optimizations for modern GPU hardware.
%We further investigated batch-oriented strategies and the integration of PDHG-based components into a feasibility-pump algorithm, see, e.g.,~\cite{berthold2019ten}.
The resulting algorithm is neither a straightforward GPU port of an existing heuristic nor a monolithic GPU-only method; rather, it is a coordinated heterogeneous portfolio with components designed to benefit from one another.
 
Our computational study is organized around this interaction between methodology and hardware.
We quantify the effect of moving from a CPU-only configuration to progressively richer CPU--GPU portfolios and compare \chap against leading GPU-based solvers and commercial CPU-based MIP solvers.

The remainder of the paper is organized as follows.
%Section~\ref{sec:TBD} gives an overview of the competition setting.
Section~\ref{sec:basics} introduces the main ingredients underlying our approach.
Section~\ref{sec:contribution} describes the individual heuristic components of our framework, their interaction in detail, and highlights our contributions.
Section~\ref{sec:experiments} reports computational results on the competition benchmark set.

% \section{The 2026 MIP Competition}
% \label{sec:overview}

\section{Prior Work}
\label{sec:basics}

Prior work on GPU-accelerated primal heuristics includes representative, often problem-specific approaches such as population-based metaheuristics (evolutionary algorithms~\cite{Soca2010,Pinel2013,Abdelatti2020}, ant colony optimization~\cite{Weiss2011,Skinderowicz2016,Uchida2012,Fingler2014}, particle swarm optimization~\cite{Zan2014,Dali2015}), local search~\cite{Schulz2013,VanLuong2011} and tabu search~\cite{Janiak2008,Czapiski2011,Bukata2013,Szymon2013,Souza2017,Abdelkafi2019}, among others; see also the survey papers~\cite{Schulz2013,Krmer2013,Boyer2013,Tan2016,Boyer2017,Essaid2018,Cheng2019,BernalNeira2024}. More general heuristic implementations for MIP include~\cite{Abbas2021,corduk25GpuAcceleratedPrimalHeuristics,Wei2025,kempke2026FP}.

\paragraph{Tabu Search}%\label{sec:tabu}

A central component of our heuristic framework is a GPU-native \emph{tabu search} designed to efficiently explore large neighborhoods in parallel. 
Tabu search is a well-established heuristic for combinatorial optimization that guides local search using tabu lists to avoid cycling and to encourage diversification \cite{Glover1989TabuSearchI,Glover1990TabuSearchII}. 
A notable recent tabu-search-based heuristic for MIP is LocalMIP \cite{Lin25LocalMip}.
%Tabu search algorithms generally follow a simple scheme. 
Given a current iterate, tabu search algorithms evaluate a set of candidate moves that modify one or multiple variable values. According to a scoring function, the best admissible move is selected and applied to the current iterate. The process is then repeated until a termination criterion is met.

\paragraph{FOM-guided Fix-and-Propagate-and-Repair}
% \label{sec:fap}

\emph{Fix-and-Propagate}~\cite{Achterberg07Thesis} is a primal heuristic that iteratively selects an integer variable, fixes it to a value within its domain, and performs domain propagation reducing the domains of the remaining unfixed variables. This process continues until all integer variables are fixed. Remaining continuous variables are then assigned by solving an LP. If the LP is feasible, the result is an integer-feasible solution; otherwise, the assignment is infeasible.
To recover from infeasible fixings, we employ the \emph{Fix-and-Propagate-and-Repair} (FPR) scheme~\cite{Salvagnin2024FPR}, which supports backtracking and repair mechanisms.
To fix variables, we use LP-based strategies as described in \cite{kempke2026FP}.

% \AGil{The following describes what we do and should be part of the next section; suggestion: have a section "FOM-guided fix-and-propagate" on kempke2026FP before FPR}

% To select variables and assign values, we use the LP solution provided by a \emph{first-order-method} (FOM)~\cite{kempke2026FP}. 
% Variable selection is performed using the following strategies, all of which rely on the LP solution:

% \begin{enumerate}
%     \item FRAC: Sorts variables by fractionality, prioritizing those closest to integer values.
%     \item DUALS: Sorts variables according to the absolute dual values of the constraints in which they appear.
%     \item REDCOSTS: Sorts variables by reduced cost, prioritizing the largest values.
%     \item TYPE: Sorts variables by type, considering binary variables first, followed by general integers.
% \end{enumerate}

% For value assignment, we use RANDOM\_LP, which rounds the LP relaxation value of a variable to the nearest integer with probability proportional to its fractional distance to each bound. For example, a variable with LP value $2.7$ is rounded to $3$ with probability $0.7$ and to $2$ with probability $0.3$.

\paragraph{Feasibility Pump}
%\label{sec:fp}

The \emph{feasibility pump} (FeasPump) is a well-known primal heuristic for MIP~\cite{berthold2019ten} that alternates between solving LP relaxations and rounding solutions to obtain integer assignments. 
Originally introduced by \cite{Fischetti2005FeasibilityPump}, the method iteratively attempts to reduce the distance between fractional LP solutions and integer points until a feasible solution is found. Our implementation is based on~\cite{Fischetti2009FeasibilityPump20}, which incorporates objective function
information into the projection step to improve solution quality, and uses FPR instead of simple rounding.

% Starting from an LP relaxation solution, the FeasPump rounds the fractional values of the integer variables to obtain an integer vector.
% This integer assignment is typically infeasible with respect to the linear constraints.
% Therefore, it is projected back onto the feasible region by solving a new LP that minimizes the distance to this rounded solution while satisfying the original constraints.
% The resulting fractional solution is then rounded again, and the process repeats

% To avoid cycling and to promote diversification, perturbation mechanisms are commonly applied when the algorithm stalls. 
% These perturbations slightly modify the rounded solution before resolving the LP, allowing the search to escape local cycles.
\section{Contribution}
\label{sec:contribution}

% \AHil{description of the figure and the overall workflow}

% \AGil{The GPU part is currently a very small part of the overall workflow.  Visually, it should take at least half, ideally much more of the space, otherwise we create the impression that our approach is mostly CPU-based.}
% \AHil{This is the problem that we have and we should discuss. The majority of the stuff that we have done on the GPU is not new. cuOpt already has a TS implementation; this is only really about engineering and (small) decision choices. What is actually the most "innovative" about our approach is that we interact with the GPU and support each other. This really helps to achieve better solutions.}

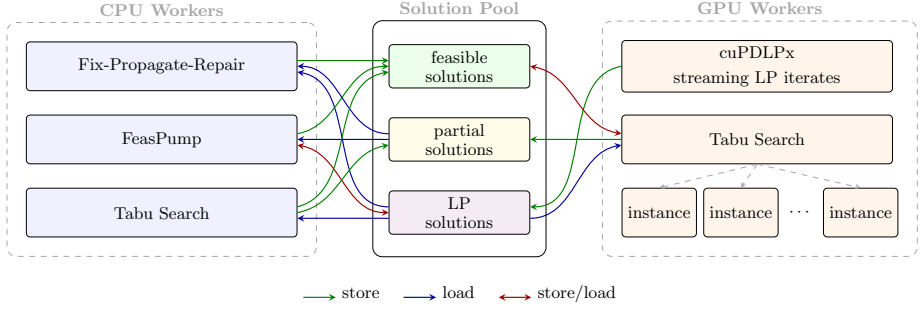
\begin{figure}[t]
\centering
\resizebox{\textwidth}{!}{%
 \begin{tikzpicture}[
  box/.style={draw, rectangle, rounded corners=2pt, minimum height=0.9cm, align=center, font=\small},
  cpubox/.style={box, fill=blue!6, minimum width=5cm},
  gpubox/.style={box, fill=orange!8, minimum width=5cm},
  instbox/.style={box, fill=orange!8, minimum width=1.1cm, font=\footnotesize},
  groupbox/.style={draw, dashed, rounded corners=5pt, inner sep=10pt, gray!70},
  store/.style={->, >=stealth, semithick, green!50!black},
  load/.style={->, >=stealth, semithick, blue!60!black},
  both/.style={<->, >=stealth, semithick, red!60!black},
]

% ============================================================
% CPU Workers (Left)
% ============================================================

\node[cpubox] (fpr) at (-5.5, 1.2) {Fix-Propagate-Repair};
\node[cpubox] (fp) at (-5.5, -0.15) {FeasPump};
\node[cpubox] (ctabu) at (-5.5, -1.5) {Tabu Search};

\node[groupbox, fit=(fpr)(fp)(ctabu),
      label={[font=\small\bfseries, gray!70]above:CPU Workers}] (cpugroup) {};

% ============================================================
% GPU Workers (Right)
% ============================================================

\node[gpubox] (pdlp) at (5.5, 1.2) {%
  cuPDLPx\\[2pt]
  {\footnotesize streaming LP iterates}%
};

\node[gpubox] (gtabu) at (5.5, -0.15) {Tabu Search};

% Instance boxes below Tabu Search
\node[instbox] (c2) at (5.2, -1.5) {instance};
\node[instbox, left=0.12cm of c2] (c1) {instance};
\node[font=\small, right=0.08cm of c2] (cdots) {$\cdots$};
\node[instbox, right=0.08cm of cdots] (cn) {instance};

\node[groupbox, fit=(pdlp)(gtabu)(c1)(cn),
      label={[font=\small\bfseries, gray!70]above:GPU Workers}] (gpugroup) {};

% ============================================================
% Solution Pool (Center)
% ============================================================

% Pool compartments at fixed positions
\node[draw, rectangle, rounded corners=2pt, minimum width=2.6cm, minimum height=0.8cm,
     fill=green!8, align=center, font=\footnotesize]
(poolfeas) at (0, 1.2) {feasible\\[-1pt] solutions};

\node[draw, rectangle, rounded corners=2pt, minimum width=2.6cm, minimum height=0.8cm,
     fill=yellow!10, align=center, font=\footnotesize]
(poolpart) at (0, -0.15) {partial\\[-1pt] solutions};

\node[draw, rectangle, rounded corners=2pt, minimum width=2.6cm, minimum height=0.8cm,
     fill=violet!8, align=center, font=\footnotesize]
(poollp) at (0, -1.5) {LP\\[-1pt] solutions};

% Outer pool frame
\path let
  \p1=(cpugroup.north), \p2=(cpugroup.south),
  \p3=(gpugroup.north), \p4=(gpugroup.south),
  \n{top}={max(\y1,\y3)}, \n{bot}={min(\y2,\y4)},
  \n{mid}={0.5*\n{top}+0.5*\n{bot}},
  \n{h}={\n{top}-\n{bot}}
in
  node[draw, rectangle, rounded corners=5pt, inner sep=8pt,
       minimum width=3.2cm, minimum height=\n{h},
       anchor=center,
       label={[font=\small\bfseries, gray!70]above:Solution Pool}]
  (pool) at (0,\n{mid}) {};

% ============================================================
% Tabu Search -> instances (dashed)
% ============================================================

\draw[->, >=stealth, dashed, gray!70] (gtabu.south) -- (c1.north);
\draw[->, >=stealth, dashed, gray!70] (gtabu.south) -- (c2.north);
\draw[->, >=stealth, dashed, gray!70] (gtabu.south) -- (cn.north);

% ============================================================
% LEFT SIDE: CPU Workers <-> Pool compartments
% ============================================================

% --- FPR ---
% FPR stores to feasible
\draw[store] ([yshift=3pt]fpr.east) -- ([yshift=3pt]poolfeas.west);
% FPR loads from partial (route down)
\draw[load] ([yshift=3pt]poolpart.west) to[out=180,in=0] ([yshift=0pt]fpr.east);
% FPR loads from LP (route down)
\draw[load] ([yshift=3pt]poollp.west) to[out=180,in=-10] ([yshift=-3pt]fpr.east);

% --- FeasPump ---
% FeasPump stores to feasible (route up)
\draw[store] ([yshift=3pt]fp.east) to[out=10,in=180] [yshift=0pt](poolfeas.west);
% FeasPump loads from partial
\draw[load] ([yshift=0pt]poolpart.west) -- ([yshift=0pt]fp.east);
% FeasPump bidirectional with LP
\draw[both] ([yshift=-3pt]fp.east) to[out=-10,in=180] ([yshift=0pt]poollp.west);

% --- CPU Tabu ---
% CPU Tabu stores to feasible (route up)
\draw[store] ([yshift=3pt]ctabu.east) to[out=15,in=190] ([yshift=-3pt]poolfeas.west);
% CPU Tabu stores to partial (route up)
\draw[store] ([yshift=0pt]ctabu.east) to[out=10,in=190] ([yshift=-3pt]poolpart.west);
% CPU Tabu loads from LP
\draw[load] ([yshift=-3pt]poollp.west) -- ([yshift=-3pt]ctabu.east);

% ============================================================
% RIGHT SIDE: GPU Workers <-> Pool compartments
% ============================================================

% --- cuPDLPx ---
% cuPDLPx stores to LP (route down)
\draw[store] ([yshift=0pt]pdlp.west) to[out=190,in=0] ([yshift=3pt]poollp.east);

% --- GPU Tabu ---
% GPU Tabu stores to feasible (route up)
\draw[both] ([yshift=3pt]gtabu.west) to[out=170,in=-10] ([yshift=0pt]poolfeas.east);
% GPU Tabu stores to partial
\draw[store] ([yshift=0pt]gtabu.west) -- ([yshift=0pt]poolpart.east);
% GPU Tabu loads from LP (route down)
\draw[load] ([yshift=-3pt]poollp.east) to[out=0,in=190] ([yshift=-3pt]gtabu.west);

% ============================================================
% Legend
% ============================================================

\node[anchor=north, font=\footnotesize] at ([yshift=-0.4cm]pool.south) {%
  \tikz{\draw[->, >=stealth, semithick, green!50!black] (0,0) -- (0.6,0);} store \quad
  \tikz{\draw[->, >=stealth, semithick, blue!60!black] (0,0) -- (0.6,0);} load \quad
  \tikz{\draw[<->, >=stealth, semithick, red!60!black] (0,0) -- (0.6,0);} store/load%
};

\end{tikzpicture}
}
\caption{Architecture of \chap.
CPU workers (left) and GPU workers (right) communicate through a shared solution pool (center).
The solution pool stores feasible solutions, partial (infeasible) solutions, and LP solutions.
On the GPU, cuPDLPx solves the LP relaxation at increasing iteration checkpoints, while the tabu search runs multiple independent instances in parallel.
}
\label{fig:architecture}
\end{figure}

\Cref{fig:architecture} depicts \chap, which combines CPU-based and GPU-based workers that interact through a shared solution pool. On the CPU side, several heuristic components run in parallel across multiple threads, enabling concurrent exploration of the search space. On the GPU side, cuPDLPx~\cite{lu2025cupdlpx} is employed to solve the LP relaxation; once this phase is completed, a GPU-based tabu search is initiated. 
All components communicate via the solution pool, facilitating the exchange of information and thereby supporting both diversification and intensification during the search.
% \AGil{For the same reason, I would put the GPU section first and ideally split it into two.  Then the warmstart section, then solution pool and summary/repetition of overall orchestration.}
% \GMil{I moved the GPU part first, followed by a small section about streaming LP iterates}
\subsection{GPU Tabu Search}

Tabu search on GPU has a long history~\cite{Janiak2008,Czapiski2011,Bukata2013,Szymon2013,Souza2017,Abdelkafi2019,corduk25GpuAcceleratedPrimalHeuristics}, among others. The most notable recent addition for MIP being~\cite{corduk25GpuAcceleratedPrimalHeuristics}, implementing LocalMIP~\cite{Lin25LocalMip} on GPU. Our tabu search implementation is based on the framework described in~\cite{Lin25LocalMip}. The following section focuses specifically on our modifications and the key distinctions from that baseline. %Aspects not covered here follow \cite{Lin25LocalMip} directly. %We denote $[n] := \{1,\dots, n\}$

Our MIPs have the form $\min_{l \leq x \leq u} \{c^T x : Ax \leq b, x_i\in\mathbb{Z}\ \forall i \in \mathcal{I} \}$ with $x, l, u, c \in \mathbb{R}^n$, variables, objective, lower and upper variable bounds, $b \in \mathbb{R}^m$ the constraint right-hand sides, $(a_{ij}) = A\in\mathbb{R}^{m\times n}$ the constraint matrix, and $\mathcal{I}\subset \{1,\dots,n\}$ the set of integer variables.

%\NKil{The section on Scoring Functions seems too verbose to me. Imo, as we do nothing different to LocalMIP, we could simply say so and cite them here? Also, why is everything called a 1-opt something here? This seems wrong. Especially as the meaning was never defined. The introduction is also broken.}
\paragraph{Scoring Functions} We assign weights $w_1,\dots,w_m\in\mathbb{R}$ to each constraint. Our candidate scoring mechanism follows~\cite{Lin25LocalMip}: given an incumbent solution $\bar{x}\in\mathbb{R}^n$ and a proposed move of variable $j\in\{1,\dots,n\}$ from $\bar{x}_j$ to $\hat{x}_j$, for each constraint $i\in\{1,\dots,m\}$ with $\sum_{j=1}^n a_{ij} x_j \leq b_i$ let
\begin{align*}
    \bar{y}_i := \sum_{k=1}^n a_{ik}\bar{x}_k \quad\text{and}\quad y_{ij} := \sum_{k\in\{1,\dots, n\}\setminus\{j\}} a_{ik}\bar{x}_k + a_{ij}\hat{x}_j.
\end{align*}
The penalty function for constraint $i$ and a move changing variable $j$ is
\begin{align*}
    p_{ij}(\hat{x}_j,\bar{x}):=
    \begin{cases}
        -w_i & \text{if } \bar{y}_i \leq b_i \text{ and } y_{ij} > b_i,\\
        w_i & \text{if } \bar{y}_i > b_i \text{ and } y_{ij} \leq b_i,\\
        0.5\,w_i & \text{if } \bar{y}_i > b_i,\; y_{ij} > b_i,\; \text{and } y_{ij} < \bar{y}_i,\\
        -0.5\,w_i & \text{if } \bar{y}_i > b_i,\; y_{ij} > b_i,\; \text{and } y_{ij} > \bar{y}_i,\\    
        0, & \text{otherwise.}
    \end{cases}
\end{align*} 
The interpretation of the scoring function: you get a reward or penalty equal to the constraint's weight when a constraint transitions between satisfied and unsatisfied (in either direction), and a reward or penalty equal to half the weight when an unsatisfied constraint becomes more or less satisfied.The total score of the move is given by
% \begin{align*}
   $ s_j(\hat{x}_j, \bar{x}) := \sum_{i=1}^m p_{ij}(\hat{x}_j,\bar{x})$.
% \end{align*}

\paragraph{Best Shift Moves} 

In contrast to \cite{Lin25LocalMip}, we employ only one type of move, the \textit{best shift move}. Given an index $j$ and a current iterate $\bar{x} \in \mathbb{R}^n$, let $\mathcal{D}_j:= \{x_j |\ l_j \leq x_j \leq u_j \text{ and } x_j \in \mathbb{Z} \text{ if } j \in \mathcal{I}\}$. A best shift move of $x_j$ at incumbent $\bar{x}$ is a change of $\bar{x}_j$ to any
\begin{equation}\label{eq:best_move}
    \hat{x}_j \in \arg\max_{x_j \in \mathcal{D}_j} s_j(x_j, \bar{x}).
\end{equation}
We remark that within this framework, it is possible that $\hat{x}_j = \bar{x}_j$. For binary variables, the only move is a flip of the binary value. For non-binary variables, we follow \cite{Luteberget2023FeasJump}: given a constraint $\sum_{j=1}^n a_{ij} x_j \leq b_i$ (w.l.o.g. $a_{ij} > 0$), we define the \emph{breakpoint} $t_{ij}$ implied on $x_j$ by constraint $i$ as
$$
\frac{b_i - \sum_{k \neq j} a_{ik}\bar{x}_k}{a_{ij}} =: t_{ij}.
$$
For $j\in\mathcal{I}$, we replace $t_{ij}$ with $\lfloor t_{ij} \rfloor$. Following a similar argument to \cite{Luteberget2023FeasJump}, the score $s_j(., \bar{x})$ is a step function in $\hat{x}_j$ with steps only at $\{t_{ij} |\ a_{ij} \neq 0, l_j \leq t_{ij} \leq u_j\} \cup \{l_j, u_j, \bar{x}_j\}$.

Prior local-search implementations consider computing (\ref{eq:best_move}) prohibitively expensive. In~\cite{Luteberget2023FeasJump} the jump value is recomputed for the most recently moved variable, enforcing $\hat{x}_j \neq \bar{x}_j$. LocalMIP~\cite{Lin25LocalMip} evaluates shifts against breakpoints from a sampled subset of violated and satisfied constraints. CuOpt~\cite{corduk25GpuAcceleratedPrimalHeuristics} evaluates all binary flips and breakpoints from violated constraints, along with a sampled subset of satisfied ones. In contrast, our implementation always computes the global optimal move for each variable w.r.t. (\ref{eq:best_move}) using the fact that $s_j(\cdot, \bar{x})$ needs to be evaluated at discrete points only. Algorithm~\ref{alg:bestshift} describes the computation of the best shift move for a single non-binary variable $x_j$ using a sort, scan, and reduce scheme (well-known GPU primitives that admit massively parallel execution~\cite{kirk2016programming}). The algorithm is motivated by the line sweep algorithm in \cite{Luteberget2023FeasJump} and consists of four phases.

\begin{algorithm}[]
\caption{Best Shift via Sort-Scan-Argmax}\label{alg:bestshift}
\begin{algorithmic}[1]
\Statex \textbf{Input:} Column $j$ with incumbent $\bar{x}_j$, bounds $[l_j, u_j]$, constraints $\sum_k a_{ik}x_k \leq b_i$ with weights $w_i$
\Statex \textbf{Output:} Best shift $\hat{x}_j$ and score $s_j$
\State $\beta \gets 0$; \; $\alpha \gets 0$; \; $B \gets \{(l_j, -1, 0),\; (u_j, -1, 0)\}$
\For{each constraint $i$ with $a_{ij} \neq 0$}
    \State $t_{ij} \gets (b_i - \sum_{k \neq j} a_{ik}\bar{x}_k) / a_{ij}$
    \If{$x_j$ integer-constrained} $t_{ij} \gets \lfloor t_{ij} \rfloor$ if $a_{ij} > 0$, \; $\lceil t_{ij} \rceil$ if $a_{ij} < 0$ \EndIf
    \If{$a_{ij} < 0$} \Comment{Imposes $x_j \geq t_{ij}$ for feasibility}
        \If{$\bar{x}_j < t_{ij}$} $\beta \mathrel{{-}{=}} \tfrac{1}{2}w_i$; $\alpha \mathrel{{+}{=}} w_i$; $B = B \cup \{(t_{ij}, -1, +\tfrac{1}{2}w_i)\}$
        \ElsIf{$\bar{x}_j > t_{ij}$} $\beta \mathrel{{-}{=}} w_i$; $B = B \cup \{(t_{ij}, -1, +w_i)\}$
        \Else{} $\beta \mathrel{{-}{=}} w_i$; $\alpha \mathrel{{+}{=}} w_i$
        \EndIf
    \Else \Comment{Imposes $x_j \leq t_{ij}$ for feasibility}
        \If{$\bar{x}_j > t_{ij}$} $\beta \mathrel{{+}{=}} w_i$; $\alpha \mathrel{{-}{=}} w_i$; $B = B \cup \{(t_{ij}, -1, 0), (t_{ij}, +1, -\tfrac{1}{2}w_i)\}$
        \ElsIf{$\bar{x}_j < t_{ij}$} $B = B \cup \{(t_{ij}, -1, 0)$, $(t_{ij}, +1, -w_i)\}$
        \Else{} $\alpha \mathrel{{-}{=}} w_i$ 
        \EndIf
    \EndIf
\EndFor
\State Sort $B$ lexicographically. Let $B^1, \dots, B^\ell$ be the ordered set of triplets.
\State $P_i \gets \sum_{k=1}^iB^k_3$
\State $\sigma_i \gets \beta + P_i + \alpha \cdot [B^k_1 > \bar{x}_j]$ for all $i = 1,\dots,\ell$.
\State $\hat{x}_j, s_j \gets \arg\max\{\sigma_i : B^i_1 \in [l_j, u_j], i = 1,\dots,\ell\}$
\end{algorithmic}
\end{algorithm}

In the first phase (lines~2--12), breakpoints are generated from each constraint. Each breakpoint $t_{ij}$ is the value of $\hat{x}_j$ at which constraint $i$ is tight. Two scalars accumulate across all constraints: $\beta$, the total reward as $\hat{x}_j \to -\infty$, and $\alpha$, the score adjustment at the incumbent. The adjustment arises from the asymmetry in the reward function: for a violated constraint, moving left from $\bar{x}_j$ worsens the violation (scoring $-\tfrac{1}{2}w_i$) while moving right improves it (scoring $+\tfrac{1}{2}w_i$). The offset $\alpha$ captures this sign flip and is applied to entries with $t_{ij} > \bar{x}_j$. We mark the $t_{ij}$  with $-1$ since the reward is credited at the breakpoint. When a breakpoint coincides with the incumbent ($t_{ij} = \bar{x}_j$), no entry is emitted; only the score change is recorded through $\beta$ and $\alpha$.

For constraints with $a_{ij} > 0$ (imposing $x_j \leq t_{ij}$), the score is unaffected up to $t_{ij}$ and drops after it. To handle this discontinuity, two entries are emitted at position $t_{ij}$: the first with second entry $-1$ and zero score change, ensuring that the score at $t_{ij}$ itself remains correct, and the second with a marker $+1$ carrying the actual penalty that should be applied after we pass $t_{ij}$. 

Phase 2 of the algorithm is a lexicographic sort over $B$. The sort order in line~13 guarantees that for a given $t_{ij}$ entries with $B^i_2 = -1$ precede ones with $B^j_2 = +1$. This guarantees that the score at $t_{ij}$ is evaluated in line~14 before the penalty for $x_j > t_{ij}$ is applied.

In the third phase is an inclusive scan which accumulates the score deltas. Each entry's final score is computed as $\beta + P_i + \alpha \cdot [\,t_i > \bar{x}_j\,]$.% where $\beta$ and $\alpha$ are broadcast to all entries. 

In the fourth phase, the argmax over entries within the variable bounds and excluding the incumbent yields the best shift move $\hat{x}_j$ and its score via a parallel reduction. A positive score indicates an improving move; if no improvement exists, the best shift is the incumbent.

\paragraph{Memory Considerations} Binary variables constitute the majority of variables in the public test set for the competition, making efficient evaluation of their flip scores critical. This has also been observed in~\cite{corduk25GpuAcceleratedPrimalHeuristics}, where, to optimize performance, flip scores are lazily reevaluated, recomputing only variables affected by the most recent move. We compute scores from scratch at each iteration. However, we implement four alternative engineering enhancements.

First, the penalty $p_{ij}$ depends on the constraint activity $\bar{y}_i$ and the right-hand side $b_i$ only through the residual $r_i = \bar{y}_i - b_i$. Storing $r_i$ directly reduces two global memory reads per nonzero to one.

Second, we normalize all constraints to the form $\sum_{j=1}^n a_{ij} x_j \leq b_i$ during preprocessing. This eliminates branching over constraint types in the evaluation kernel. On GPUs, divergent branches within a warp force serialization of both paths, making uniform control flow essential for performance.

Third, evaluating flip scores requires iterating over the nonzero entries of the constraint matrix. This iteration can be organized either row-wise (CSR order) or column-wise (CSC order), each with a distinct memory-access trade-off. In the row-wise approach, threads in a warp process nonzeros from the same row, sharing the same residual and constraint weights - these values can be broadcast or read from L1 cache without generating redundant global memory traffic. The cost is a scattered read of the incumbent and an atomic addition to a per-variable score array. In the column-wise approach, threads in a warp share the same variable and scores can be accumulated via warp-level reduction without atomics. However, each thread now references a different row, requiring scattered reads of both the residual (a double) and the constraint weight (a float). Our experiments confirm that the row-wise approach is faster when we need to recompute the score of all binary columns. A consequence of this design choice is that selective reevaluation of scores for only affected variables, as proposed in~\cite{corduk25GpuAcceleratedPrimalHeuristics}, becomes expensive in CSR order since affected columns are not contiguous.

Fourth, we store the binary incumbent as a bitset and pin it to the L2 cache via \texttt{cudaAccessPropertyPersisting}. Since our solver runs multiple parallel climbers sharing the same device, L2 capacity is at a premium. The bitset representation reduces the per-climber footprint, ensuring that all climbers' incumbents remain L2-resident even on GPUs with limited L2 (approximately 50\,MB on the H100). This further supports the row-wise design choice, as the column-wise approach would require pinning the much larger residual and weight arrays to L2, which is impractical given the limited cache capacity.

\paragraph{Kernel Specialization} The implementation of Algorithm~\ref{alg:bestshift} relies on performant sort, scan, and reduction primitives. We make three remarks.

First, since columns in typical MIP constraint matrices are short, we implement custom kernels where multiple short vectors (of length 2, 4, 8, or 16) are packed into a single warp and sorted simultaneously using bitonic sorting networks~\cite{kipfer2005improved} via warp-level primitives~\cite{cudawarpprim}. Bitonic sorting networks have a fixed comparison topology that is independent of the input data, thereby avoiding control-flow divergence within a warp. We remark that while CUB provides segmented reduce and scan at the warp level~\cite{cub}, to the best of our knowledge no existing CUDA library natively supports sorting multiple independent segments cooperatively within a single warp.

Second, we dispatch to specialized kernels based on vector length to minimize wasted cycles and maximize memory locality. Vectors of length up to 64, 128, or 256 are packed into 256-thread blocks, where each thread handles one element of its assignend vector and communication occurs through shared memory. Vectors of length up to 1024 are assigned to a single block of the nearest enclosing power-of-two size. Vectors longer than 1024 are sorted using grid-wide primitives spanning multiple blocks. A native alternative, launching one block per short vector without packing, was observed in preliminary testing to be slower.

Last, instead of dispatching all our kernels individually, we use CUDA Gra\-phs~\cite{cudaGuide} to record multiple tabu search iterations and their respective kernel submissions, eliminating most of the kernel submission overhead and allowing for parallel kernel execution.

\paragraph{Parallel Tabu Search Instances}

To improve both hardware utilization and search diversification, we execute multiple tabu search instances in parallel. Each instance maintains its own current solution, tabu list, and constraint weights, and therefore explores an independent search trajectory in the solution space.

Running many search instances in parallel allows the algorithm to fully exploit the massive parallelism of modern GPUs while reducing the risk of prematurely converging to a single region of the solution space.

% The parallel structure also enables asynchronous algorithmic updates. While some search instances are evaluating moves, global components of the heuristic framework, such as restart decisions or solution-pool communication (see Section~\ref{sec:solutionpool}), can be updated without interrupting the remaining instances. Search instances may also be restarted when their progress stagnates or when the diversity between trajectories drops below a predefined threshold.

% \paragraph{Diversity scoring}

% To prevent different search instances from converging to nearly identical solutions, we maintain a diversity measure between their current iterates. For two solutions $x$ and $y$, we compute a generalized Hamming distance over the discrete variables
% \[
% d(x,y) = \sum_{i \in \mathcal{D}} \mathbf{1}[x_i \neq y_i],
% \]
% where $\mathcal{D}$ denotes the set of binary and integer variables. 
% This measure counts the number of discrete-variable assignments that differ between two solutions.

% The diversity score of a search instance is computed relative to the other active instances, and instances whose solutions become too similar to others get restarted. This mechanism encourages the parallel tabu searches to explore different regions of the solution space.

\subsection{Solution Pool}\label{sec:solutionpool}
% nk: I know technically we have 2 pools - but I think describing them as one pool makes the presentation simpler. AH: agree
The different components of the heuristic framework communicate through a shared solution pool. By centralizing solution management in this way, the pool acts as a lightweight coordination mechanism that allows the different components of the framework to cooperate while still operating largely independently. As a result, locking between CPU threads is minimal and restricted primarily to writing solutions.
As shown in \Cref{fig:architecture}, the solution pool is divided into three parts: 
LP solutions generated by cuPDLPx, and FeasPump (see \Cref{sec:lpstreaming}), feasible incumbent solutions, and infeasible solutions used to guide FPR and tabu search (see \Cref{sec:repair}).
Each heuristic has access to the current best-known objective value, which it uses to impose a cutoff constraint that guides the search toward improving solutions.

% The pool stores both feasible incumbents and promising infeasible solutions used to guide FPR and FeasPump.
% The solution pool serves several purposes. 
% First, it provides a global objective cutoff shared by all algorithmic components. 
% Both TS and FPR are using the best objective to add this additional cutoff constraint to avoid finding feasible but not improving solutions.
% Second, the solution pool enables the exchange of information between the tabu search and the FPR and FeasPump framework. 
% Tabu search runs may be restarted from solutions sampled from the pool, thus allowing promising regions of the search space to be explored further.  
% This enables tabu search to act as an \emph{improvement heuristic}. 
% Conversely, infeasible solutions with good objective values discovered during tabu search are also stored and can be used as starting reference points for the FPR and FeasPump procedures.
% Both, which attempts to repair these solutions and drive them towards feasibility.  This allows FPR to act as a \emph{repair heuristic}.

\subsection{Streaming LP Iterates}\label{sec:lpstreaming}

% A key design choice is that the LP relaxation solved by cuPDLPx on the GPU does not block the CPU workers.

Since the convergence of the Primal-Dual Hybrid Gradient (PDHG) method (the base of cuPDLPx) can be relatively slow, relying on full convergence before using the LP solution would introduce unnecessary latency.
Moreover, because the LP solution is primarily used to guide primal heuristics such as FPR, FeasPump, and Tabu Search, it is not clear that a highly accurate LP solution is required (and \cite{mexi2023scylla,kempke2026FP} indicate it is not).
Following the idea of concurrent processing discussed by Rothberg for crossover \cite{rothberg2025concurrent}, we adopt a streaming approach: 
Instead of waiting for convergence, we stream intermediate LP iterates at checkpoints after $10^2$, $10^3$, $10^4$, and $10^5$ PDHG iterations, warm-starting each phase from the previous iterate.
Each checkpoint produces a snapshot consisting of the primal solution, dual values, and reduced costs, which is immediately deposited into the solution pool.
These intermediate solutions are used to initialize the heuristics with essentially no latency.

% \subsection{FPR}
% However, the FPR workers start simultaneously with the PDLP solver computing the LP relaxation, so the LP solution is initially unavailable. 
% In the meantime, we use the LP-free variable selection and assignment strategies proposed by~\cite{Salvagnin2024FPR}.
% %
% Additionally, for TSP-like instances, we use the additional LP-free strategy TSP that ranks variables by their appearance in eq-cliques\cite{hansknecht2021dynamic}. \AH{@NK this is the correct Hansknecht paper?}

% We introduce the assignment strategy \texttt{RANDOM_GUIDED}. 
% Similar to its LP equivalent, variables are fixed based on probabilities derived from their (fractional) values in the TS solution. 
% For selection, we employ the \texttt{RowViolation} and \texttt{Type} strategies.

\subsection{Using a FeasPump-FPR Combination as a Repair Heuristic}
\label{sec:repair}

Both FeasPump and FPR access partial solutions and intermediate LP solutions from the solution pool in order to derive feasible incumbent solutions.

\paragraph{FeasPump} The FeasPump worker operates in two phases. In the first phase, FeasPump evaluates the intermediate LP solution from \Cref{sec:lpstreaming} and attempts to generate feasible solutions. If feasibility cannot be achieved, the resulting solutions are nevertheless added to the partial solutions pool to serve as guidance for the FPR workers.
In the second phase, for each partial (infeasible) solution submitted by the tabu search (TS), a limited number of FeasPump iterations is performed to move the solution toward feasibility.

\paragraph{FPR} Besides LP-based strategies applied to the (intermediate) LP solution and TSP-based strategies~\cite{hansknecht2021dynamic}, we use FPR as a repair mechanism for solutions generated by FeasPump or the tabu search. Similar to the LP-based strategies in \cite{kempke2026FP}, we fix the values of selected variables accordingly and, if necessary, apply repairing steps to restore feasibility.

In summary, the repair mechanism exploits infeasible solutions produced by the tabu search, which may first be moved closer to feasibility by FeasPump. This process diversifies the FPR approach and strengthens the ability of our heuristics to escape local minima.
\section{Experiments}
\label{sec:experiments}

% \AHil{@Gioni: we discovered some divergence when calculating the PI with Nils and Gen. We fixed it. But your results deviate from Nils/Gen's. We should check that all scripts are consistent.}

\paragraph{Setup}
We evaluate our approach on the 50 instances provided by the competition organizers, presumably chosen s.t. they are challenging for current CPU- and GPU-based solvers, with a time limit of 5 minutes. All CPU components run on an Intel(R) Xeon(R) Platinum 8481C processor, while the GPU components (cuPDLPx~\cite{lu2025cupdlpx,lu2024restarted} and GPU tabu search) run on an NVIDIA H100 HBM3 GPU. In our experiments, we used Xpress 9.8.0, Gurobi 12.0.1 and cuOpt 25.10. Xpress is used by CHAP for solving LP relaxations after FPR has fixed all integer variables.

%\paragraph{GPU Utilization} We assessed the GPU utilization of our using the profiling tools provided by NVIDA~\cite{Nsystems,Ncompute}. Figure~\ref{fig:profiling} shows the Nsight Systems profiler output for three launched graph instances of our tabu search kernels, running on instance\_49, one of the larger test instances. Both occupancy and active streaming multiprocessors are high. Our kernels are mainly shared memory bound, caused by the sort operation within Algorithm~\ref{alg:bestshift}, leading to a relatively low instruction throughput. We note that occupancy and throughput degrade when considering smaller instances.

\paragraph{CHAP Components Evaluation}  To assess the contribution of the GPU components, we compare three configurations: a pure CPU variant, a variant that adds the GPU-accelerated LP solve via cuPDLPx, and the full portfolio that additionally includes the GPU tabu search.
\begin{table}[t]
    \centering
    \begin{tabular}{l@{\hskip 2em}r@{\hskip 2em}r@{\hskip 2em}r@{\hskip 2em}r}
      \toprule
      Configuration & Found & Wins & Gap\% & PI \\
      \midrule
      Pure CPU        & 47          &  3 & 19.57\% & 72.59 \\
      CPU + GPU LP    & 46          &  13 & 12.17\% & 51.01 \\
      CPU + GPU LP + GPU Tabu & 47  & \textbf{22} & \textbf{8.84\%} & \textbf{41.84} \\
      \bottomrule
    \end{tabular}
    \caption{Comparison of three configurations on 50 instances with a 5-minute time limit. Wins are based on best final solution quality, counted out of 47 instances solved by at least one configuration (9 ties). Gap\% (shifted geometric mean, shift 1.0) and \emph{Primal Integral} (Definition 2.3~\cite{Berthold2013}) are computed with respect to the best solution found by Gurobi or Xpress within 8\,h.}
    \label{tab:gpu_ablation}
\end{table}
%
%\GMil{Also the long table in the appendix is a bit confusing. In table 1 we compute gap and PI wrt the best solution of Gurobi or Xpress in 8h, In table 2 we compare to VBS(5min Gurobi and Xpress), and in the appendix compare to 5min Gurobi, and VBS (Gurobi or Xpress in 8h). Suggestion: In the appendix we name it Best Known 8h and omit the 5min Gurobi columns}
Table~\ref{tab:gpu_ablation} reports a comparative analysis of three configurations evaluated on competition instances with a time limit of 5 minutes per instance. 
The columns report the number of instances solved (“Found”), the number of wins, as well as the shifted geometric means of the primal gap (Gap\%) and primal integral (PI). 
As references for Gap and PI, we compare against the best solution computed by either Gurobi and Xpress within an 8 hour time limit. The instance-wise best solution value of these runs can be found in the Appendix~\ref{sec:appendix:overall_results}, Table~\ref{tab:appendix} in the column ``VBS 8 h''.

All combinations find a near-identical amount of solutions up to one instance, which generally proved unstable in our experiments. Replacing the CPU-based PDHG with its GPU-based counterpart a substantial improvement in solution quality. In particular, the average gap is reduced from 19.57\% to 12.17\%, and the PI decreases from 72.59 to 51.01, indicating faster convergence toward high-quality solutions.
A more pronounced improvement is observed when GPU-accelerated tabu search is additionally incorporated (CPU + GPU LP + GPU Tabu). 
This configuration achieves a significantly higher number of wins (22), clearly outperforming the other approaches. 
This is reflected in the lowest gap (8.84\%) and PI (41.84), indicating both superior final solution quality and improved convergence behavior over time.

% \GMil{Should we omit table 2? We can just mention it in the text.}
% \begin{table}[h!]
%       \centering
%       \begin{tabular}{l@{\hskip 2em}r@{\hskip 2em}r@{\hskip 2em}r@{\hskip 2em}r}
%         \toprule
%         & Found & Wins & Only  \\
%         \midrule
%             Ours         & 47 & 10 & 2  \\%             Virtual best & 46 & 33 & 1 \\
%             \bottomrule
%       \end{tabular}
%       \caption{Solution quality comparison of our best configuration against the
%   virtual best of Gurobi, Xpress with heuristic emphasis, and cuOpt (default and
%   heuristics mode) on 50 instances with a 5-minute time limit. Wins count instances
%   where one solver found a strictly better objective; Only counts instances
%   solved by only that solver. (two ties)}
%       \label{tab:vs_virtualbest}
% \end{table}

\begin{table}[t]
      \centering
      \begin{tabular}{l@{\hskip 2em}r@{\hskip 2em}r@{\hskip 2em}r}
        \toprule
        & Found & Wins & Only \\
        \midrule
        \chap            & 47 & 11 & 1 \\
        Gurobi           & 44 & 25 & 1 \\
        cuOpt (heur)     & 43 &  6 & 0 \\
        \bottomrule
      \end{tabular}
      \caption{Solution quality comparison of the best \chap configuration against
  Gurobi and cuOpt (\texttt{--heuristics-only} mode) on 50 instances with a 5-minute time limit. Wins count
  instances where a solver found the unique best objective; Only counts
  instances solved exclusively by that solver. (6~ties, 2~unsolved by all).}
      \label{tab:vs_commercial}
    \end{table}
%\TBil{How would the Wins statistic look like when we compared Wins/Losses against VB and cuOpt individually? So one number X:Y for us against the commercial solvers and one for us against cuOpt?}
%\GMil{Against cuopt 24:19 (4 ties), against VB 17:30 (1 tie)}
%
\paragraph{Solver Comparison} Table \ref{tab:vs_commercial} compares our best configuration against the commercial CPU-based MIP solver Gurobi and against the GPU-based MIP solver Nvidia cuOpt in \texttt{--heuristics-only} mode under a 5-minute time limit. 
The results are reported in terms of the number of solved instances (``Found''), the number of wins (instances with a uniquely best objective value), and the number of instances solved exclusively by a given approach (``Only'').

Our method achieves the highest number of solved instances (47), outperforming both the virtual best commercial solver (45) and cuOpt (43). 
In addition, it is able to find one solution where all other solvers fail to find one.
In terms of wins, however, the virtual best attains the highest count (25), reflecting the complementary strengths of the underlying commercial solvers in achieving superior objective values. 
Our method secures 11 wins, demonstrating competitive solution quality.
In a direct comparison against cuOpt, our framework finds a better solution 24 times, cuOpt 19 times, with 5 ties, demonstrating that the two approaches complement each other nicely.

% \begin{table}[h!]
% \centering
% \begin{tabular}{l|rrrrr}
% \hline
% Solver & Solved & Raw SGM & Raw Mean & Pen SGM & Pen Mean \\
% \hline
% cuopt\_default & 39 & 7.56 & 18816.49 & 13.73 & 14698.87 \\
% cuopt\_heuristics & 42 & 19.93 & 14533908.40 & 25.92 & 12208499.06 \\
% mini\_pump\_polish & 50 & 8.42 & 38.43 & 8.42 & 38.43 \\
% \hline
% \end{tabular}
% \caption{Summary Statistics (based on 50 instances with a reference value). Penalized gaps assign 100\% to unsolved instances. Raw statistics consider only solved instances.}
% \end{table}

% Summary Statistics (based on 50 instances with a reference value)

% Solver                 Solved      Raw SGM     Raw Mean      Pen SGM     Pen Mean
% --------------------------------------------------------------------------------
% cuopt\_default              32        11.23       140.86        25.15       126.15
% cuopt\_heuristics           38         8.51        28.02        15.77        45.30
% main                       50        13.60        52.57        13.60        52.57

\section{Conclusion}

We presented a hybrid framework \chap that combines methodological and engineering innovations to exploit CPU and GPU resources in a coordinated way for MIP heuristics.
On the methodological side, our approach integrates multiple strategies, warmstarting, portfolio ideas, and improvement mechanisms into a unified cross-platform search procedure.
On the engineering side, we presented specialized kernels, LP iterate streaming, memory pinning and efficiency-improved move evaluations that make these ideas effective in practice.

The computational results on the competition benchmark are encouraging: \chap outperforms cuOpt and commercial CPU-based solvers on the given instance set, highlighting the potential of coordinated cross-platform primal heuristics. Integrating the approach into a state-of-the-art MIP solver is a promising future research direction.

\bigskip

\noindent
\textbf{Acknowledgements}
The work for this article has been conducted in the Research Campus MODAL funded by the German Federal Ministry of Research, Technology, and Space (BMFTR) (fund numbers 05M14ZAM, 05M20ZBM,\\05M2025).

\newpage
\appendix

\section{Overall results}\label{sec:appendix:overall_results}
%\AHil{PI is currently against VBS}

% CHAP Obj ist ausgelesen aus /scratch/gcp1/nkempke/results/paper_final/default/*
% PI 
% Gurobi 5 Minuten
% https://www.mixedinteger.org/2026/competition/instance_5m_bounds.txt
% VBS is the best out of Gurobi and Xpress 
% https://git.zib.de/gtjusila/mipcc2026/-/wikis/Baseline-Runs/Gurobi-50GB-RAM,-8-Hours,-8-Threads
% https://git.zib.de/gtjusila/mipcc2026/-/wikis/Baseline-Runs/XPRESS-50GB-RAM,-8-Hours,-8-Threads
\begin{table}[H]
\centering
\small
\resizebox{0.9\textwidth}{!}{%   % width = \textwidth, height = auto (preserve aspect ratio)
\begin{tabular}{lrrrrrr}
\toprule
& \multicolumn{2}{c}{\chap}& \multicolumn{2}{c}{Gurobi 5m} & \multicolumn{2}{c}{VBS 8 h} \\
\cmidrule(r){2-3}
\cmidrule(r){4-5}
\cmidrule(l){6-7}
Instance & Obj & PI & Obj & Bound & Obj & Bound \\
\midrule
01 & -211940.97 & 11.158 & -218764.89 & -221054.62 & -218764.89 & -218767.38 \\
02 & -- & 300.000 & -357538.58 & -358714.57 & -357544.31 & -357579.87 \\
03 & 310.60 & 63.839 & 248.50 & 248.50 & 248.50 & 248.50 \\
04 & 927.60 & 12.416 & 1089.50 & 491.30 & 568.70 & 491.30 \\
05 & 971.70 & 14.583 & 1036.40 & 442.80 & 589.40 & 442.80 \\
06 & 53114.00 & 14.734 & 55056.00 & 38161.00 & 46616.00 & 38527.61 \\
07 & 80825.00 & 15.713 & 80985.00 & 75477.00 & 76812.00 & 75506.96 \\
08 & 44955.00 & 35.359 & 40056.00 & 38848.00 & 39122.00 & 38995.00 \\
09 & 15.00 & 106.388 & 10.00 & 10.00 & 10.00 & 10.00 \\
10 & 11.00 & 36.555 & 10.00 & 7.00 & 10.00 & 10.00 \\
11 & 561.00 & 25.190 & 729.00 & 148.00 & 422.00 & 389.60 \\
12 & 214.00 & 34.220 & 193.00 & 193.00 & 193.00 & 193.00 \\
13 & 4205.00 & 28.285 & 5053.00 & 14.00 & 2482.00 & 1571.00 \\
14 & 693643.00 & 102.194 & 495373.00 & 484942.00 & 490492.00 & 487858.00 \\
15 & 1016442.00 & 113.125 & 678131.00 & 669447.54 & 674011.00 & 672812.01 \\
16 & 46 & 221.51 & -- & 5.00 & 44.00 & 6.00 \\
17 & -- & 300.000 & -- & 1.00 & -- & 3.00 \\
18 & 17.00 & 3.592 & 16.00 & 16.00 & 16.00 & 16.00 \\
19 & 6.00 & 51.845 & 5.00 & 5.00 & 5.00 & 5.00 \\
20 & 10.00 & 35.224 & 7.00 & 4.00 & 7.00 & 5.00 \\
21 & 11.00 & 13.893 & 10.00 & 7.00 & 10.00 & 7.00 \\
22 & 211.01 & 70.793 & 171.50 & 159.74 & 169.44 & 166.33 \\
23 & 352.48 & 214.741 & 100.84 & 100.84 & 100.84 & 100.84 \\
24 & 436.51 & 11.233 & 431.24 & 425.10 & 431.24 & 431.20 \\
25 & -32.00 & 0.138 & -32.00 & -32.00 & -32.00 & -32.00 \\
26 & -90.00 & 1.440 & -86.00 & -110.00 & -86.00 & -108.00 \\
27 & -145.00 & 1.927 & -92.00 & -324.00 & -145.00 & -148.00 \\
28 & 269141.51 & 31.236 & 373739.50 & 45621.01 & 208447.50 & 56165.82 \\
29 & 931303.78 & 46.902 & -- & 17801.67 & -- & 113510.53 \\
30 & -- & 300.000 & -- & 0.000 & --- & 12443.64 \\
31 & -2487.11 & 10.650 & -2507.62 & -2514.56 & -2505.87 & -2509.70 \\
32 & -3443.12 & 27.947 & -3493.74 & -3503.17 & -3494.89 & -3496.96 \\
33 & -2260.03 & 36.050 & -331.70 & -10679.23 & -10536.74 & -10675.07 \\
34 & 118.00 & 59.851 & 100.00 & 100.00 & 100.00 & 100.00 \\
35 & 172.00 & 50.078 & 148.00 & 139.00 & 148.00 & 148.00 \\
36 & 414.00 & 59.854 & 328.00 & 270.00 & 320.00 & 287.00 \\
37 & 3.00 & 150.075 & 2.00 & 0.00 & 2.00 & 2.00 \\
38 & 245.00 & 172.626 & 89.00 & 0.00 & 42.00 & 0.00 \\
39 & 490.00 & 17.067 & 464.00 & 439.00 & 464.00 & 464.00 \\
40 & 591.00 & 13.022 & -- & 461.00 & 524.00 & 524.00 \\
41 & 776.00 & 14.851 & -- & 503.00 & 634.00 & 516.00 \\
42 & 24956015.12 & 6.192 & 24942375.27 & 24941304.54 & 24942293.49 & 24940450.33 \\
43 & 24371169.11 & 13.061 & 24253497.15 & 24203429.80 & 24252006.89 & 24251105.07 \\
44 & 57080573.56 & 77.487 & 57045735.99 & 56965375.53 & 57026496.19 & 57022946.64 \\
45 & 4.00 & 114.653 & 4.00 & 4.00 & 4.00 & 4.00 \\
46 & 4.00 & 137.013 & 56.00 & 1.00 & 1.00 & 1.00 \\
47 & 146.00 & 249.049 & 116.00 & 1.00 & 22.00 & 1.00 \\
48 & 153219161.81 & 6.425 & 152287062.76 & 152272618.81 & 152287056.67 & 152286169.89 \\
49 & 127362905.84 & 8.329 & 269331393.81 & 123721906.89 & 126091275.99 & 125528411.68 \\
50 & 130427528.91 & 8.885 & 130086220.74 & 124660126.38 & 129906242.85 & 129329721.70 \\
\bottomrule
\end{tabular}%
}
\caption{\footnotesize Comparison of objective values obtained by \chap with Gurobi (5-minute time limit) and the virtual best solver (VBS), combining Gurobi and Xpress (8-hour time limit). PI denotes the primal gap, using the best solution found by either the VBS or \chap as reference.}
\label{tab:appendix}
\end{table}

% \section{GPU Profiling}
% \begin{figure}[htbp]
%   \centering
%   \includegraphics[width=\linewidth]{figures/profiling1.png}
%     \caption{Nsight Systems timeline showing GPU activity, SM occupancy, and stream utilization during a Tabu Search iteration on an NVIDIA H100.}

%   \label{fig:profiling}
% \end{figure}
\newpage
% \newpage
% \begin{itemize}
%   \item Intro (AG, TB, SP)
%   \begin{itemize}
%     \item Considerations: Verteilungsproblemen, Autotuning, andere Hardware-Überlegungen (TK)
%   \end{itemize}

%   \item Contributions:
%   \begin{itemize}
%     \item Multiple Climbers: Scores (NK)
%     \item Portfolio-Approaches mit noch nicht getesteten Kombinationen (AH)
%     \begin{itemize}
%       \item Solution Pool (AH)
%     \end{itemize}
%     \item Warmstarting (AH)
%     \item LocalMIP als Improvement-Heuristik: FPR + Local-MIP (NK)
%     \item FJ-Sampling in LocalMIP; andere Neuigkeiten in LocalMIP (NK)
%     \item Drei specialized Kernel: ganz lang, ein Block, spezialisiert (1 mehr) (GT)
%     \item Move-to-move operation: Können vollständig auswerten, wo vorher heuristisches Sampling nötig war (GT)
%     \item Engineering: pinning (GT)
%     \item Implementierte Strategien: Batch-FP (GM)
%     \item PDHG in FPR/FP (GM)
%   \end{itemize}

%   \item Visualisierung / Diagramm: Interaktion der Komponenten (NK + AH)

%   \item Experimente:
%   \begin{itemize}
%     \item High-Level Experiment: CPU-classic -- CPU-PDHG -- GPU-PDHG -- Portfolio (GM)
%     \item Wir vs.\ Nvidia (GT)
%     \item LocalMIP -- original vs.\ unser LocalMIP (GT)
%   \end{itemize}
% \end{itemize}
\bibliographystyle{splncs04}
\bibliography{bibliography}

\end{document}